# A Heuristic Approach to Two Level Boolean Minimization Derived From Karnaugh Mapping

**Ethan L. Childerhose**
Iroquois Ridge H.S.
Oakville, Ontario
info@ethanchilderhose.com

**Jingzhou Liu**
Iroquois Ridge H.S.
Oakville, Ontario
jason_liu1116@hotmail.com

June 19, 2019

## Abstract

The following paper presents a heuristic method by which sum-of-product Boolean expressions can be simplified with a specific focus on the removal of redundant and selective prime implicants. Existing methods, such as the Karnaugh map and the Quine-McCluskey method [1, 2], fail to scale since they increase exponentially in complexity as the quantity of literals increases, doing as such to ensure the solution is algorithmically obtained. By employing a heuristic model, nearly all expressions can be simplified at an overall reduction in computational complexity. This new method was derived from the fundamental Boolean laws, Karnaugh mapping, as well as truth tables.

## 1 Introduction

Given an unsimplified Boolean expression in the form of sum-of-products, there exist implicants that are both redundant and non-trivial. For the purpose of this paper, a trivial solution is considered to be a Boolean expression that can be simplified using the basic laws of Boolean algebra without introducing extraneous product terms to the original expression. An example of a trivial solution is as follows. Given two terms $\bar{A}\bar{B} + A\bar{B}$, this expression can be simplified to $\bar{B}$ with the use of the OR distributive law and the complement law as seen in Table 1.1.

| Step | Rule | Expression |
|------|------|------------|
| 1 | Given | $\bar{A}\bar{B} + A\bar{B}$ |
| 2 | Distributive | $\bar{B}(\bar{A} + A)$ |
| 3 | Complement | $\bar{B}$ |

**Table 1.1:** Trivial Absorption Proof

However, there exists a set of Boolean expressions that cannot be efficiently solved by following fundamental Boolean rules. An example of one such expression is $\bar{A}\bar{C}+\bar{A}B+BC$. According to the truth table demonstrated in Table 1.2, $\bar{A}\bar{C}+\bar{A}B+BC$ is logically identical to $\bar{A}\bar{C}+BC$ in spite of having no trivial solution. When this expression is represented on a Karnaugh map as shown in Figure 1.1, the redundant term, $\bar{A}B$, becomes apparent.

| A | B | C | $\bar{A}\bar{C}+\bar{A}B+BC$ | $\bar{A}\bar{C}+BC$ |
|---|---|---|---|---|
| 0 | 0 | 0 | 1 | 1 |
| 0 | 0 | 1 | 0 | 0 |
| 0 | 1 | 0 | 1 | 1 |
| 0 | 1 | 1 | 1 | 1 |
| 1 | 0 | 0 | 0 | 0 |
| 1 | 0 | 1 | 0 | 0 |
| 1 | 1 | 0 | 0 | 0 |
| 1 | 1 | 1 | 1 | 1 |

**Table 1.2:** Truth Table of $\bar{A}\bar{C}+\bar{A}B+BC$ and $\bar{A}\bar{C}+BC$

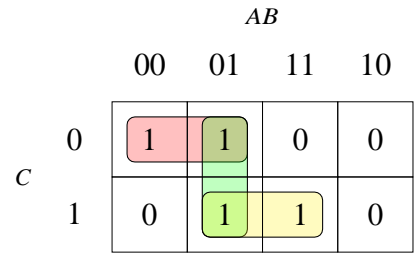

**Figure 1.1:** K-map of $\bar{A}\bar{C}+\bar{A}B+BC$

While there exist several techniques that can accomplish the removal of logic redundancies, there are certain inadequacies that pertain to each method. For instance, one can achieve logic minimization algebraically by adding extraneous implicants that aid with the merging of terms. However, in the example demonstrated in Table 1.3 and Table 1.4, in which multiple selective terms can be expunged, determining which selective implicant leads to the most optimal expression is unclear. As a result, all of the possible paths have to be carried out in order to determine the optimum expression, which becomes progressively in-feasible as the number of implicants increases.

| Step | Rule | Expression |
|---|---|---|
| 1 | Given | $\bar{A}\bar{C}\bar{D}+\bar{A}B\bar{C}+ACD+ABD+B\bar{C}D$ |
| 2 | Identity | $\bar{A}\bar{C}\bar{D}+\bar{A}B\bar{C}+ACD+ABD+B\bar{C}D \cdot 1$ |
| 3 | Complement | $\bar{A}\bar{C}\bar{D}+\bar{A}B\bar{C}+ACD+ABD+B\bar{C}D \cdot (A+\bar{A})$ |
| 4 | Distributive | $\bar{A}\bar{C}\bar{D}+\bar{A}B\bar{C}+ACD+ABD+AB\bar{C}D+\bar{A}B\bar{C}D$ |
| 5 | Absorption | $\bar{A}\bar{C}\bar{D}+\bar{A}B\bar{C}+ACD+ABD$ |

**Table 1.3:** Non-Optimal Algebraic Boolean Solution

| Step | Rule | Expression |
|------|------|------------|
| 1 | Given | $\bar{A}\bar{C}\bar{D} + \bar{A}B\bar{C} + ACD + ABD + B\bar{C}D$ |
| 2 | Identity | $\bar{A}\bar{C}\bar{D} + \bar{A}B\bar{C} \cdot 1 + ACD + ABD \cdot 1 + B\bar{C}D$ |
| 3 | Inverse | $\bar{A}\bar{C}\bar{D} + \bar{A}B\bar{C} \cdot (D + \bar{D}) + ACD + ABD \cdot (C + \bar{C}) + B\bar{C}D$ |
| 4 | Distributive | $\bar{A}\bar{C}\bar{D} + \bar{A}B\bar{C}D + \bar{A}B\bar{C}\bar{D} + ACD + ABCD + AB\bar{C}D + B\bar{C}D$ |
| 5 | Absorption | $\bar{A}\bar{C}\bar{D} + ACD + B\bar{C}D$ |

**Table 1.4:** Optimal Algebraic Boolean Solution

Other techniques, such as Karnaugh maps and the Quine–McCluskey algorithm [1, 2], suffer a similar issue as the algebraic approach since the computational cost of conducting those methods grows exponentially with the number of variables. The purpose of this paper is to demonstrate a new heuristic and expandable approach that accomplishes two level Boolean simplification in a more efficient manner than current methods.

# 2 The Tail-Eliminate Method

## 2.1 Procedure

1. Expand the Boolean expression to ensure every possible prime implicant is included.
2. Generate the Tail-Eliminate map from the expanded expression
3. Check if any of the following end conditions are satisfied. If one is met, then the most optimal expression has been achieved; otherwise, continue with the procedure.
   (a) Tail Quotient of all implicants is greater than 0
   (b) Tail Quotient of all implicants equal
4. Determine the prime implicant(s) with the lowest tail quotient. These implicants are selective prime implicants.
5. Determine the prime implicant(s) with the highest tail quotient. These are essential prime implicants.
6. Verify if any of the selective prime implicants overlap with a tail implicant (an essential prime implicant).
   (a) If none of the implicants overlap with any of the selected prime implicants, then the simplified expression has been reached.
   (b) If yes, remove that implicant from the expression.
7. Repeat back to step 2

## 2.2 Expansion

The first step to conducting two-level Boolean minimization is to ensure that all prime implicants are included in the Boolean expression. Any expansion method can be used to generate prime implicants given a truth table or an expression. For the research conducted in this paper, expansion methods in BOOM [3], were employed to generate a set of prime implicants.

## 2.3 Determining Overlaps Between Boolean Terms

Overlap is key to determining logic redundancies. Overlap is most simply defined as the area on a Karnaugh map where two prime implicants exhibit crossover. This is demonstrated in Figure 2.1, the Karnaugh mapping of the expression $\bar{A}\bar{C} + \bar{A}B$.

$$AB$$

|   | 00 | 01 | 11 | 10 |
|---|----|----|----|----|
| C=0 | 1 | 1 | 0 | 0 |
| C=1 | 0 | 1 | 0 | 0 |

**Figure 2.1:** Overlap Between $\bar{A}\bar{C}$ and $\bar{A}B$

On this Karnaugh map, there is a total of one overlap between the two implicants. To understand the way in which the number of overlaps can be mathematically calculated, a more rigorous explanation for overlaps is required.

An overlap, more acutely, is a partial redundancy within the truth table of an expression. A partial redundancy occurs when one given input produces a TRUE output for two different product terms. Visualizing partial redundancies is best demonstrated through the Boolean expression $\bar{A}\bar{C} + \bar{A}B$ when represented in a truth table as shown in Table 2.1.

| A | B | C | $\bar{A}\bar{C}$ | $\bar{A}B$ | $\bar{A}\bar{C} + \bar{A}B$ |
|---|---|---|------------------|------------|------------------------------|
| 0 | 0 | 0 | 1 | 0 | 1 |
| 0 | 0 | 1 | 0 | 0 | 0 |
| **0** | **1** | **0** | **1** | **1** | **1** |
| 0 | 1 | 1 | 0 | 1 | 1 |
| 1 | 0 | 0 | 0 | 0 | 0 |
| 1 | 0 | 1 | 0 | 0 | 0 |
| 1 | 1 | 0 | 0 | 0 | 0 |
| 1 | 1 | 1 | 0 | 0 | 0 |

**Table 2.1:** Truth Table of $\bar{A}\bar{C} + \bar{A}B$

In the case of $\bar{A}\bar{C} + \bar{A}B$, both terms yield a TRUE when A is FALSE, B is TRUE, and C is FALSE, hence causing a partial redundancy. Given that specific input, only one product term is required to return a TRUE since the two terms in sum of product form undergo the OR operation.

From here, the basis for the way in which overlap is determined between two terms becomes easier to grasp. As seen in the truth table, when the logical AND operation is conducted between the two product terms $\bar{A}\bar{C}$ and $\bar{A}B$, the resulting product term $\bar{A}B\bar{C}$ only returns TRUE when both inputs from the two product terms are TRUE. Consequently, the AND operation creates a resultant product term ($\bar{A}B\bar{C}$) that implies the overlap area of the two product terms as shown in Figure 2.2. Therefore, the number of TRUE outputs of the resultant product term is equal to the number of overlaps, or partial redundancies.

| A | B | C | $\bar{A}\bar{C}$ | $\bar{A}B$ | $\bar{A}\bar{C} \cdot \bar{A}B = \bar{A}B\bar{C}$ |
|---|---|---|---|---|---|
| 0 | 0 | 0 | 1 | 0 | 0 |
| 0 | 0 | 1 | 0 | 0 | 0 |
| **0** | **1** | **0** | **1** | **1** | **1** |
| 0 | 1 | 1 | 0 | 1 | 0 |
| 1 | 0 | 0 | 0 | 0 | 0 |
| 1 | 0 | 1 | 0 | 0 | 0 |
| 1 | 1 | 0 | 0 | 0 | 0 |
| 1 | 1 | 1 | 0 | 0 | 0 |

**Table 2.2:** Truth Table of $\bar{A}\bar{C} \cdot \bar{A}B$

|    | AB 00 | 01 | 11 | 10 |
|----|----|----|----|----|
| C 0 | 1 | 1 | 0 | 0 |
| 1 | 0 | 1 | 0 | 0 |

**Figure 2.2:** K-Map of the Resultant Product Term $\bar{A}B\bar{C}$

In order to determine the number of TRUE outputs of the resultant product term, a function, $o(n, a) = 2^{n-a}$, can be defined, with $n$ being the total quantity of unique literals in the original Boolean truth table and $a$ being the quantity of unique literals in the resultant product term. Since the expression, $\bar{A}\bar{C}+\bar{A}B$, has three unique literals (A, B, and C) and the resultant product term, $\bar{A}B\bar{C}$, also contains three unique literals, the function $o(n, a)$ outputs 1. Thus, there is one overlap between the implicant $\bar{A}\bar{C}$ and the implicant $\bar{A}B$.

However, it is worth noting that the overlap function, $o(n, a)$, should not be used when the variable, $a$, equals zero, which only occurs when one product term includes a complimented version of a literal that exists in another product term. For instance, the implicants $\bar{A}B$ and $AC$ both share the literal A, except for the fact that A is complimented in only one out of the two implicants. When applying the AND operation between $\bar{A}B$ and $AC$, it is realized that the resultant product term is $A\bar{A}BC$, which can be further simplified to 0 according to the Law of Complement. In this case when there is no literal in the resultant product term, the number of overlaps between the two product terms is zero.

The number of overlaps is crucial to determining logic redundancies because when all of the TRUE outputs of a term are overlapped by other prime implicants, the term is considered to be a redundant prime implicant and can therefore be expunged from the Boolean expression. By compiling and examining the number of overlaps of each prime implicant through the use of a Tail-Eliminate map, prime implicants can be differentiated into essential prime implicants and redundant prime implicants.

## 2.4 Generating the Tail-Eliminate Map

The Tail-Eliminate map is a table that represents and calculates the accumulated number of overlaps between each term in a sum-of-products Boolean expression. The horizontal axis and the vertical axis are labeled as the product terms of an expression that contains redundant prime implicants and essential prime implicants. The following example, $\bar{A}\bar{C}\bar{D}+\bar{A}B\bar{C}+B\bar{C}D+ABD+ACD$, demonstrates the way in which the rest of the Tail-Eliminate map should be set up.

|  | AB |  |  |  |
|---|---|---|---|---|
| CD | 00 | 01 | 11 | 10 |
| 00 | 1 | 1 | 0 | 0 |
| 01 | 0 | 1 | 1 | 0 |
| 11 | 0 | 0 | 1 | 1 |
| 10 | 0 | 0 | 0 | 0 |

**Figure 2.3:** K-Map of $\bar{A}\bar{C}\bar{D}+\bar{A}B\bar{C}+B\bar{C}D+ABD+ACD$

| x | $\bar{A}\bar{C}\bar{D}$ | $\bar{A}B\bar{C}$ | $B\bar{C}D$ | $ABD$ | $ACD$ |
|---|---|---|---|---|---|
| $\bar{A}\bar{C}\bar{D}$ | x | 1 | 0 | 0 | 0 |
| $\bar{A}B\bar{C}$ | 1 | x | 1 | 0 | 0 |
| $B\bar{C}D$ | 0 | 1 | x | 1 | 0 |
| $ABD$ | 0 | 0 | 1 | x | 1 |
| $ACD$ | 0 | 0 | 0 | 1 | x |
| Total Overlaps | 1 | 2 | 2 | 2 | 1 |

**Table 2.3:** TE-Map of $\bar{A}\bar{C}\bar{D}+\bar{A}B\bar{C}+B\bar{C}D+ABD+ACD$

The highlighted value in Table 2.3 is calculated through the overlap function $o(n, a)$ between $\bar{A}B\bar{C}$ and $\bar{A}\bar{C}\bar{D}$, which are the horizontal term and the vertical term of that cell respectively. The remaining values in the table are computed in a similar manner; each column is then summed to produce the total quantity of overlaps for that column's implicant. If the total number of overlaps for any particular implicant equals zero, that implicant needs to be omitted from the Tail-Eliminate map.

The total number of overlaps is used to determine the tail implicant through the calculations of the Tail Quotient, both of which are concepts elaborated in the following section.

## 2.5 Tail Implicant and Tail Quotient

For the purpose of this research paper, a tail implicant is defined to be a subset of essential prime implicants with the highest tail quotient across all other prime implicants in a Boolean expression. The tail quotient of an implicant is calculated using the equation, $t = h - o_t$, with $h$ being the quantity of TRUE outputs implied by the implicant and $o_t$ being the total number of overlaps of the implicant. In order to determine the tail implicant, a new row that displays the tail quotient of each implicant can be added beneath the Total Overlaps row in the Tail-Eliminate map, as shown in Table 2.4. In Table 2.4, the implicant $\bar{A}\bar{C}\bar{D}$ as well as the implicant $ACD$ hold the largest tail quotient and are thus considered to be the tail implicant.

| x | $\bar{A}\bar{C}\bar{D}$ | $\bar{A}B\bar{C}$ | $B\bar{C}D$ | $ABD$ | $ACD$ |
|---|---|---|---|---|---|
| $\bar{A}\bar{C}\bar{D}$ | x | 1 | 0 | 0 | 0 |
| $\bar{A}B\bar{C}$ | 1 | x | 1 | 0 | 0 |
| $B\bar{C}D$ | 0 | 1 | x | 1 | 0 |
| $ABD$ | 0 | 0 | 1 | x | 1 |
| $ACD$ | 0 | 0 | 0 | 1 | x |
| Total Overlaps | 1 | 2 | 2 | 2 | 1 |
| Tail Quotient | 1 | 0 | 0 | 0 | 1 |

**Table 2.4:** Tail-Eliminate Mapping of $\bar{A}\bar{C}\bar{D} + \bar{A}B\bar{C} + B\bar{C}D + ABD + ACD$

The implicants with the lowest tail quotient are defined to be a selective prime implicant; they are prime implicants that have the potential to be removed, given that the procedure's end conditions (will be discussed in the next section) have not yet been met. To determine whether or not one particular selective prime implicant can be expunged, the concept of tail implicants must be applied. More specifically, a selective prime implicant is redundant and can be removed if it overlaps with the tail implicant in any given step of iteration, otherwise, the selective prime implicant should remain untouched since only one implicant is removed per iteration. In the example shown above, implicants $\bar{A}B\bar{C}$, $B\bar{C}D$, and $ABD$ are all selective prime implicants; however, only $\bar{A}B\bar{C}$ and $ABD$ can be removed since they overlap with the tail implicant $\bar{A}\bar{C}\bar{D}$ or $ACD$ whereas $B\bar{C}D$ does not. In this situation where two of the selective prime implicants can be removed, one of the two is selected and expunged. Through the employment of tail terms, the correct order of steps can be found avoiding the incorrect solution as shown in Table 1.3.

## 2.6 End Conditions

As redundant terms are determined and removed, it is necessary to constantly reassess the expression to verify that the most optimal solution has not been reached. There exists two end conditions that allow for the determination of the most optimal solution. If all tail quotients are equal or all tail quotients are above 0, then the most optimal solution has been reached and all implicants are now essential prime. The basis behind these conditions is best understood through a Tail-Eliminate map.

| x | $\bar{A}\bar{C}\bar{D}$ | $B\bar{C}D$ | $ACD$ |
|---|---|---|---|
| $\bar{A}\bar{C}\bar{D}$ | x | 0 | 0 |
| $B\bar{C}D$ | 0 | x | 0 |
| $ACD$ | 0 | 0 | x |
| Total Overlaps | 0 | 0 | 0 |
| Tail Quotient | 2 | 2 | 2 |

**Table 2.5:** Tail-Eliminate Mapping of $\bar{A}\bar{C}\bar{D} + B\bar{C}D + ACD$

Table 2.5 is the final map when simplifying the expression $\bar{A}\bar{C}\bar{D} + \bar{A}B\bar{C} + B\bar{C}D + ABD + ACD$. After the two terms $\bar{A}B\bar{C}$ and $ABD$ are removed, the two end conditions can be clearly seen. The first condition that is met is that all tail quotients are above 0. When any given implicants tail quotient is above 0, that implicant's quantity of overlaps does not exceed the minterms within that implicant and therefore there will be minterms that are not overlapped, making the implicant a prime essential implicant. When this occurs across all implicants within an expression, that expression must be fully simplified. There exist situations, however, where the tail quotient of implicants remains at 0 or below. Given this situation, the solution is found when the tail quotients are equal. When tail quotients are equal, each implicant is equally likely to be a redundant term and as such the expression is considered to be optimized.

## 3 Conclusion

The Tail-Eliminate method is currently a heuristic method that approximates the most optimal solution as it has not yet been proven to simplify all possible expressions. We hope to determine the degree of accuracy of the Tail-Eliminate method and develop it into an algorithm if possible. We hope to further test the Tail-Eliminate method and benchmark it against other existing heuristic methods. Both of these goals will be achieved by creating programs capable of automatically running expressions through the Tail-Eliminate method. Up to date programs for the Tail-Eliminate method can be found at https://github.com/EthChil/BooleanSimplification.